\documentclass[12pt]{amsart}

\usepackage{amssymb, amsfonts}

\theoremstyle{plain}

\newcommand{\beq}{\begin{equation}}
\newcommand{\eeq}{\end{equation}}

\newtheorem{theor}{Theorem}[section]

\newtheorem{rem}{Remark}[section]

\numberwithin{equation}{section}

\begin{document}

\begin{center}
{\large {\bf Valiron-Titchmarsh Theorem for Positive Temperatures} }\\

\vspace{.5cm}

{\bf J. B. Lacay }
\vspace{.25cm}

Bronx Community College of the City University of New York
Mathematics and Computer Science Department \\
2155 University Avenue, Bronx, NY 10453, USA \vspace{.1cm}\\
\medskip

\end{center}

\vspace{.1cm}

{\bf {Abstract}.} In this note, we prove an analog of the Valiron-Titchmarsh theorem for positive temperatures. 
\vspace{.25cm}

\begin{center} {\bf {Key words}:} Entire functions, Counting function, Parabolic cylinder. \vspace{.25cm}
\end{center}


\section{Introduction and Statement of Results}

\emph{Let $f(z)$ be an entire function of order $\rho <1$ with only negative zeros.}\

\[\emph{$\log f(r)\approx \frac{\pi}{\sin \pi \rho}r^{\rho}$ as $r\rightarrow \infty.$ }\]

\noindent If the counting function $n(r)$ of the zeros of $f$ satisfies $n(r)\approx r^{\rho}, r\rightarrow \infty$, then, it is immediate that

\[\emph{$ \log f(r)\approx \frac{\pi}{\sin \pi \rho}r^{\rho}, r\rightarrow \infty.$}\]

\noindent The converse statement of Tauberian nature was proved independently by Valiron and Titchmarsh \cite{Tit}, and it is often referred as the Valiron-Titchmarsh theorem on entire functions with negative zeros; its current exposition can be found,  in  \cite [Lec 12-13]{Lev2}. For the history and generalizations of the Valiron-Titchmarsh theorem (see, e.g., \cite{Lev1} or \cite{Khe1} and the references therein. In \cite{Str}, it was shown that, for an entire function of non- integer order with zeros on the negative real half-line, the existence of the asymptotics of a certain form for one of the functions along some ray implies the existence of certain asymptotics for the counting function of the zeros. In \cite{Khe1}, the author showed that for an entire function $f(z)$ of finite order $\rho$ with all its zeros lie on a finite collection of rays in the interior of a sector $S: \alpha<arg z<\beta$, such a function is called completely regular along $arg z=\theta (\alpha ,\beta )$. If the $ {\lim_{r \rightarrow \infty }} r^{- \rho}  \log |f(r e^{i \theta}| $ exists, then $f(z)$ is of completely regular growth on the whole plane. Recently, \cite{Khe2} extended the Drasin complement to the Valiron-Titchmarsh theorem and showed that if $u$ is a subharmonic function of this class and order $0<\rho <1 $, then the existence is  the $ {\lim_{r \rightarrow \infty }} \log \frac {u(r)}{N(r)} $. In \cite{Khe&Lac}, the problems under consideration examined the relationship of the initial data $f$ and that of the solution $u$. From the main theorem, we proved the interesting corollary that $u(x,t)= \alpha x^{\rho} + \bar{\bar{o}}\big(r^{\alpha}\big)$ for each $t \rightarrow \infty$ if and only if $ f(y)=\alpha y^{\alpha}+\bar{\bar{o}}\big(y^{\alpha}\big), y \rightarrow \infty $. In \cite {Lac}, we considered two-term analogs of the Valiron-Titchmash theorem for the temperatures.

In this note, we prove an analog of the Valiron-Titchmarsh theorem  for positive temperatures, i.e for positive solutions of the heat equation. 

\beq \frac{\partial u(x,t)}{\partial t}=\kappa
\frac{\partial^2 u(x,t)}{\partial x^2} \eeq in the slab $\mathcal{S}^T= \mathbb{R}^2 \times (0,T) $, where the constants $\kappa>0$ and $0< T\leq \infty.$
Thereafter, these solutions are called temperatures. 

It is known \cite [p.57]{Wat} that a positive temperature in $S^T$ has the Gauss-Weierstrass representation
\beq u(x,t)= \frac{1}{4\pi \kappa t} \int_{\mathbb{R}^2} e^{{-\frac{1}{4\kappa t}}\|x-y\|^2 } d\mu(y).\eeq

where $x = (x_1, x_2)$, $x=(r \cos \theta, r \sin \theta),  r\geq 0, 0\leq \theta<2 \pi $, and $y = (y_1, y_2)$, $y=(s \cos \phi, s \sin \phi),  s\geq 0, 0\leq\phi<2 \pi $ in spherical coordinates respectfully. Here, $\|\centerdot \|$ is the Euclidean norm in $\mathbb{R}^2, $ and $d \mu$ is a non-negative function on $\mathbb{R}^2$.

It is known \cite[p.57, Theor 2.10]{Wat} that if $u$ is real-valued and continuous up to the boundary  $\mathbb{R}^2 \times [0] $, then $d\mu(y)=u(y,0)dy.$ It is also known \cite[p.8]{Wat}, that under this assumption the measures $\mu$ are absolutely continuous with respect to the Lebesgue measure in $\mathbb{R}^2$.
To derive an analog of the  Valiron- Titchmarsh theorem for the positive temperatures, we assume that the function  $d\mu(y)$ is supported on the ray $\arg y=\theta_0$. Thus, in the case under consideration, (1.2) can be represented as
\beq u(x,t)= \frac{1}{4\pi \kappa t} \int_{0}^ \infty  e^{{-\frac{1}{4t}}\|x-y\|^2 } dn (s),\eeq
where now $y=(s \cos\theta_0,s \sin\theta_0)$ and $n(s)$ denotes the number of zeros in the circle $|z| \leq r$ and $n(s)=\mu(\{|y|\leq s\}) $ is the total $\mu -$measure of the disk $|y| \leq s$.  Since we are interested in asymptotic properties of temperatures, we can assume without any loss of generality that some vicinity at the origin is $\mathbb{R}^2 $ is free of the measure $d\mu$, i.e. $n(s_0)=0$ for some $s_0>0$, in particular, $n(0)=0$.
Moreover, since $u(0,t)<\infty$, we have from (1.3)
\beq \int_{0}^{\infty}  e^{{-\frac{1}{4t}}s^2 }dn(s) < \infty .\eeq 
We also assume for the rest of the note that for any $t>0,$
\beq \lim_{s\rightarrow\infty}n(s) e^{-s^2/4t}=0.\eeq
Integrating (1.4) by parts and using (1.5), we derive the representation
\beq u(0,t)= \frac{1}{4\pi \kappa t} \int_{0}^{\infty}  s n(s) e^{{-\frac{1}{4t}}s^2 }ds .\eeq 

\noindent Suppose the measure $u$ is supported on the ray $(s, \theta_0), s>0, 0 \leq \theta_0 <2 \pi.$ Thus, in (1.3) the measure $n(s)=n((0,s])$ of the semi-inteval $(0,s]$ for all $s>0.$ Integrating (1.3) by parts, we have

\beq u(x,t) =\frac {1}{2t} \int_{0}^\infty (s-rcos (\theta -\theta_0 ) )   e^{-\frac {r^2 -2r s cos(\theta -\theta_0 )+s^2}{4t} } n(s) ds     \eeq 

\[ =\frac {1}{2t} e^{-\frac {r^2}{4t}} \int_{0}^\infty (s-rcos (\theta -\theta_0 ) )   e^{-\frac {s^2 -2r s cos(\theta -\theta_0 )}{4t} } n(s) ds \]  

\begin{rem} Since we are interested in temperatures with power growth of the measure, when (1.5) is clearly valid, (1.5) is not an essential restriction for us.
 \end{rem}

\begin{theor} Let the temperature $u$ have a representation
\[u(x,t)=Ax^\alpha + v(x,t) \]
where
\[ \lim_{x \rightarrow\infty} x^{-\alpha}\int_{0}^{x} \upsilon (y,t)dy=0\]
is uniform in $t \in [0,t_0]$ for some $0<t_0<a$. Then
\[f(r)=Ar^{\alpha}+\bar{\bar{o}}\big(r^{\alpha}\big), r \rightarrow \infty. \]
\end{theor}
\noindent \textbf{Proof.} For a function $f$, denote $[f](x)=\frac{1}{2} \big ((f(x^+)+f(x^-) \big)$. By assumption$f(y)=0$ in some neighborhood of $y=0$, thus $[f](0)=0$ and by  \cite[p. 69, Theor. 6]{Wid},
\[ \lim_{t \rightarrow 0} \int_{0}^{r}u(y,t)dy=[f](r)\]
therefore
\[[f](r)=\lim_{t \rightarrow 0} \int_{0}^{r}f(y) dy+\lim_{t \rightarrow 0} \int_{0}^{r} \upsilon (y,t)dy \]
\[=A r^\alpha+\lim_{t \rightarrow 0} \int_{0}^{r}v(y,t)dy. \]
The limits as $ r \rightarrow +\infty$ and as $t \rightarrow 0^+$ can be interchanged due to the uniformity assumption, thus
\[\lim_{r \rightarrow \infty}r^{-\alpha} \lim_{t \rightarrow 0^+} \int_{0}^{r}v(y,t)dy =\lim_{t \rightarrow 0^+}\lim_{r \rightarrow \infty}r^{-\alpha}  \int_{0}^{r}v(y,t)dy=0.\]

\begin{theor} Under some asumptions on a temperature u

\[u(r)= \lim_{t {\rightarrow {0}^+} } \int_{|y|\leq r}   u(y,t)dy  \]

\noindent if $u(x,t) \approx A( \theta_0) |x|^{\rho }.$ Then

\beq  f|x|) = \alpha |x|^{\rho} +  \overline{\overline{o}} (x^\rho )\eeq   

\end{theor}

\noindent \textbf{Proof.} By \cite[Theor. 7.2]{Wat} 

\beq u(r)= \lim_{t {\rightarrow {0}^+} } \int_{|y|\leq r}   u(y,t)dy  \eeq

\[=    \lim_{t {\rightarrow {0}^+} } \int_{|y|\leq r}  A(\theta) B^\rho dy   + \lim_{t {\rightarrow {0}^+} } \int_{|y|\leq r}   u(y,t)dy , y=(y_1, y_2) =(B,\theta)  .  \]\

\noindent Suppose 
\[ \lim_{t {\rightarrow {0}^+} } r^{-s-2}  \int_{|y|\leq r}   V(s,\theta ,t)dy =0.  \]

\noindent If under some asumptions on a temperature $u$,

\[ u(r)= \lim_{t {\rightarrow {0}^+} } \int_{|y|\leq r}   u(y,t)dy . \]

\noindent Then, 

\[  \lim_{r=|x| {\rightarrow {\infty }} } \frac {u(x,t)}{|x|^{\rho +2}} = \frac {1}{\rho+2} \int_{2}^{2 \pi} A(\theta ) d \theta . \]

\noindent Since

\[ u(x,t)= A(t, \theta ) |x|^\rho + V(x,t), y=(y_1, y_2) =(B,\theta) .\]

\noindent It follows,

\[\int_{|y| \leq r} A(t, \theta) s^\rho dy = \int_0^{r} \int_0^{ 2 \pi} s A(t, \theta ) s^\rho d \theta ds \]\

\[ = \int_0^{ 2 \pi}  A(t, \theta ) \int_0^r s^{\rho +1} ds \]\

\[ =\frac {r^{\rho+2}}{\rho +2}  \int_0^{2 \pi} A(t, \theta ) d \theta \]

\noindent Now, we can state our result.

\begin{theor} Let $u(x,t)$ be a positive temperature given by (1.2) with the measure $d \mu $ supported at the ray, $arg z = \theta_0$. If $n(s_0)=0 $ and 
\beq n(s)=a_0s^{\alpha(s)}+n_1(s), s>s_0,\eeq
where the constants $a_0$ and $\alpha $ satisfy $a_0\geq0$, $\alpha>-1$, and the remainder
\[\lim_{s\rightarrow\infty}s^{-\alpha(s)} n_1(s)=0, \]
then

\beq u(x,t)=
\begin {cases}
\frac{a_0}{2 \sqrt{\pi t}} \big (r \cos (\theta- \theta_0) \big )^\alpha e^{-\frac{r^{2} \sin^2(\theta- \theta_0)}{4 t}},&\text{ $ \cos(\theta-\theta_0)>0$ }\\\\
\frac {a_0}{\pi} \Gamma(\frac {\alpha+1}{2}) 2^{\alpha -2} t^{\frac{\alpha -1}{2}} e^{-\frac{r^2}{4 t}}, &\text{ $ \cos(\theta-\theta_0)=0.$ }\\\\
\frac {a_0}{\pi} \Gamma(\alpha +1) 2^{\alpha -1} \frac{t^\alpha}{\big(r| \cos(\theta-\theta_0)| \big)^{\alpha +1}} e^{-\frac{r^2}{4 t}}, &\text{ $ \cos(\theta-\theta_0)<0.$ }
\end {cases} \eeq
\end{theor}
\noindent \textbf{Proof.} We write (1.3) as
\beq u(x,t)=u_0(x,t)+u_1(x,t) \eeq
\[\equiv \frac{a_0}{8\pi t^2} \int_{0}^{\infty} e^{{-\frac{1}{4t}}||x-y||^2 }  s^{\alpha(s)} ds + \frac{b}{8\pi t^2} \int_{0}^{\infty} e^{{-\frac{1}{4t}}||x-y||^2 }  n_1(s) ds. \]
We find the principal term of the asymptotic formula by estimating $u_0,$
\beq u_0(x,t)= \frac{a_0}{4\pi \kappa t} \int_{0}^{\infty} e^{{-\frac{1}{4t}}||x-y||^2 }  s^{\alpha(s)} ds\eeq
\[ = \frac{a_0}{4 \pi \kappa t} \int_{0}^{\infty}  e^{-\frac{(s^2-2rs \cos(\theta-\theta_0+r^2)}{4t}} s^{\alpha(s)} ds ,\]
\[ = \frac{a_0  e^{-\frac{r^2}{4t}}}{4 \pi \kappa t} \int_{0}^{\infty}  e^{-\frac{s^2-2rs \cos(\theta-\theta_0)}{4t}} s^{\alpha(s)} ds .\]
To simplify the integral $u_0$, set $s=\sqrt{2t}w,$ $\alpha=-\nu-1$,  while letting $z=-\frac{r \cos(\theta-\theta_0)}{\sqrt{2t}}$,
\[ u_0(x,t)= \frac{a_0}{4 \pi \kappa t}  e^{-\frac{r^2}{4t}} \big (\sqrt{2t}\big )^{\alpha(s)+1}\int_{0}^{\infty} w^{\alpha} e^{{-\frac{w^2}{2}}-zw}dw .\] 
The  above integral $u_0$ can be expressed through and estimated by making use of the parabolic cylinder functions, where $D_{\nu}(z)$ is the  \emph{Weber function} \cite[Section 8.3,(3)]{BatErd}, 
\beq D_{\nu}(z)=\frac{e^{-z^2/4}}{\Gamma(-\nu)}\int_0^{\infty} e^{-zt-{t^2/2}}t^{-\nu-1}dt,\; \Re \nu <0.\eeq 
Using (2.2) and (2.3),  $u_{0}$ can be written as
\[u_0(x,t)=\frac{a_0 \Gamma(\alpha(s)+1)}{\pi} 2^{\frac {\alpha-3}{2}}  t^{\frac {\alpha-1}{2}} e^{-\frac {r^2}{8t} \big(1+sin^2 (\theta-\theta_0)\big)} D_{-\alpha -1}(z),\] 
We consider here two cases,  $ \cos(\theta-\theta_0)>0$ and $ \cos(\theta-\theta_0)<0.$ We cut the $z-plane$ along the negative $x-axis$, $x=Re$$z$, thus $-\pi <\theta\leq \pi$ and fix the value $arg (1)=0.$ Then in the half-plane, $\theta_0 - \pi/2<\theta<\theta_0 + \pi/2$, $z=-\frac{r \cos(\theta-\theta_0)}{\sqrt{2t}}<0, arg z=\pi.$\\
Using the known asymptotic formulas\footnote{The corresponding formulas in \cite[Sect. 8.4]{BatErd} contain misprints - missing brackets.} \cite[p. 307]{BleHan} for a fixed value of $\nu$ and $\pi/4 <|\arg z|<5/4\pi,$ as $z \rightarrow \infty $,, we have
\[D_{\nu}(z)= z^{\nu}  e^{-\frac{1}{4}z^2} \bigg [ \sum_{n=0}^{N} \frac{(-1/2 \nu)_n (1/2-1/2 \nu)_n}{n! (-1/2z^2)^{n}} + \bar{\bar{o}}\big(|z^{2}|^{-N-1}\big) \bigg ]\] \[ - \frac {(2 \pi)^{1/2}}{\Gamma (-\nu)}e^{\nu \pi i} z^{-\nu -1} e^{\frac{1}{4}z^2} \bigg [ \sum_{n=0}^{N} \frac{(1/2 \nu)_n (1/2+1/2 \nu)_n}{n! (1/2 z^2)^{n}} + \bar{\bar{o}}\big(|z^{2}|^{-N-1}\big) \bigg ].\]\\
For $ \cos(\theta-\theta_0)>0$ , $  e^{\frac{1}{4}z^2}\gg e^{-\frac{1}{4}z^2}$,  $n=0$,
\[D_{\nu}(z)\approx  - \frac {(2 \pi)^{1/2}}{\Gamma (-\nu)}e^{\nu \pi i} z^{-\nu -1} e^{\frac{1}{4}z^2} \]
Thus,
\[u_0(x,t)=\frac{a_0}{\pi} 2^{\frac {\alpha(s)-3}{2}}  t^{\frac {\alpha(s)-1}{2}}   \big ( r \cos(\theta-\theta_0) \big )^{\alpha(s)} e^{-\frac {r^2}{4t} sin^2 (\theta-\theta_0)} \]
Similarly, for $ \cos(\theta-\theta_0)<0$
\[D_{\nu}(z)\approx  \frac {(2 \pi)^{1/2}}{\Gamma (-\nu)} e^{-\frac{1}{4}z^2} \]
\[u_0(x,t)=\frac{a_0}{\pi} 2^{\frac {\alpha-3}{2}}  t^{\frac {\alpha-1}{2}} \big ( r \cos(\theta-\theta_0) \big )^\alpha e^{-\frac {r^2}{4t}} \]
Now, we estimate $u_1$ by (1.6),
\[ u_1(x,t)=\frac{b}{4\pi \kappa  t} \int_{0}^{\infty} e^{{-\frac{1}{4t}}||x-y||^2 }  n_1(s) ds. \]
Since $\lim_{s\rightarrow\infty}s^{-\alpha(s)} n_1(s)=0.$
Therefore,
\[A e^{{-\frac{1}{4t}}||x-y||^2 } \rightarrow 0, s\rightarrow \infty,\]
 uniformly in $t.$\\

\end{document}